\documentclass[a4paper,reqno,11pt]{amsart}

\usepackage{fullpage}
\usepackage{microtype}
\usepackage[utf8]{inputenc}
\usepackage{mathtools,amssymb}
\usepackage{enumitem}
\usepackage{tikz}
\setlist[enumerate]{label=(\roman*),font=\normalfont}
\usepackage[colorlinks=true,linkcolor=blue,citecolor=magenta,urlcolor=magenta]{hyperref}
\usepackage{cleveref}
\usetikzlibrary{decorations.markings}
\usetikzlibrary{patterns,decorations.pathreplacing}

\newtheorem{thm}{Theorem}[section]
\newtheorem{prop}[thm]{Proposition}
\newtheorem{cor}[thm]{Corollary}
\newtheorem{lem}[thm]{Lemma}
\theoremstyle{definition}
\newtheorem{defi}[thm]{Definition}
\newtheorem{ex}[thm]{Example}
\theoremstyle{remark}

\crefname{thm}{Theorem}{Theorems}
\crefname{prop}{Proposition}{Propositions}
\crefname{cor}{Corollary}{Corollaries}
\crefname{lem}{Lemma}{Lemmas}
\crefname{defi}{Definition}{Definitions}
\crefname{ex}{Example}{Examples}
\crefname{rem}{Remark}{Remarks}
\numberwithin{equation}{section}

\newcommand{\kk}{\Bbbk}
\newcommand{\ZZ}{\mathbb Z}
\newcommand{\RR}{\mathbb R}
\newcommand{\eb}{\mathbf e}
\newcommand{\tb}{\mathbf t}
\newcommand{\one}{\mathbf 1}
\newcommand{\frakm}{\mathfrak m}
\newcommand{\Fill}{\operatorname{Fill}}
\DeclareMathOperator{\relint}{int}
\newcommand{\Hilb}{\operatorname{Hilb}}
\newcommand{\conv}{\operatorname{conv}}
\newcommand{\type}{r}

\begin{document}

\title{Almost Gorenstein edge rings of bipartite graphs and their $h$-vectors}
\author{Yuta Hatasa}
\address[Y. Hatasa]{Independent Researcher, Japan}
\email{hatasa.y.b8ad@gmail.com}
\author{Koji Matsushita}
\address[K. Matsushita]{Graduate School of Mathematical Sciences, The University of Tokyo, Komaba, Meguro-ku, Tokyo 153-8914, Japan}
\email{koji-matsushita@g.ecc.u-tokyo.ac.jp}

\begin{abstract}
    We characterize bipartite graphs whose edge rings are almost Gorenstein.
    Using this characterization, we recover known results for two families of graphs.
    In addition, we show that if the edge ring of a bipartite graph is almost Gorenstein and its $h$-polynomial has degree at least two, then the $h$-polynomial is palindromic except possibly for the linear and next-to-leading coefficients.
    Moreover, we construct two bipartite graphs whose edge rings have the same $h$-vector, although only one of the rings is almost Gorenstein.
\end{abstract}

\maketitle

\section{Introduction}

A Cohen--Macaulay graded ring is almost Gorenstein if it admits a degree-zero embedding into an appropriate shift of its canonical module whose cokernel is either zero or an Ulrich module.
Here,
an Ulrich module is a nonzero Cohen--Macaulay module whose multiplicity equals its minimal number of generators.
Almost Gorenstein graded rings were introduced by Goto, Takahashi, and Taniguchi (\cite{goto2015almost}) and have since been studied in various settings, including
graded rings arising from combinatorial objects (see, e.g., \cite{matsuoka2016uniformly,miyazaki2018hibi,miyazaki2023stable}).

By a theorem of Stanley, a Cohen--Macaulay standard graded domain is Gorenstein if and only if its $h$-vector is symmetric (\cite[Theorem~4.4]{stanley1978hilbert}).
This naturally raises the question of whether almost Gorensteinness is also determined by the $h$-vector.
Higashitani investigated the relationship between almost Gorensteinness and $h$-vectors (\cite{higashitani2016almost}) and asked whether there exist two Cohen--Macaulay standard graded domains with the same $h$-vector such that one is almost Gorenstein and the other is not (\cite[Problem~4.6]{higashitani2016almost}).

In this paper, we study almost Gorensteinness for edge rings of bipartite graphs, which form a class of normal Cohen--Macaulay domains.
For a graph $G$, its edge ring $\kk[G]$ is the $\kk$-subalgebra of $\kk[t_v:v\in V(G)]$ generated by $t_ut_v$ for $\{u,v\}\in E(G)$, with $\deg(t_ut_v)=1$.
Almost Gorenstein edge rings have been studied for the following families:
\begin{itemize}
    \item Complete multipartite graphs (\cite{higashitani2022levelness}).
    \item Odd cycles sharing a common vertex (\cite{higashitani2023h,bhaskara2026h}).
    \item Complete bipartite graphs with a matching deleted (\cite{hatasa2025pseudo}).
\end{itemize}

Our main result gives a complete characterization of bipartite
graphs whose edge rings are almost Gorenstein.
We also investigate the behavior of their $h$-vectors.
To state our main theorem, we recall some graph-theoretic
terminology (see \Cref{subsec:edge-rings} for details).
A graph is \textit{matching-covered} if it is connected, has at
least one edge, and every edge belongs to a perfect matching.
For a vertex $v$ of $G$, let $N_G(v)$ denote its neighborhood.
For a $2$-connected matching-covered bipartite graph $G$,
the filling construction introduced in
\cite{hatasa2026palindromicity} defines a set $\Fill(G)$
of missing edges whose addition yields a graph with a
Gorenstein edge ring.

By the block decomposition and the tensor-product results of Matsushita and Miyashita (\cite{matsushita2024multiplicity}), it suffices to consider $2$-connected bipartite graphs (see \Cref{prop:2-connected}).
\begin{thm}\label{thm:main}
    Let $G$ be a $2$-connected bipartite graph with bipartition $V(G)=X\sqcup Y$.
    Then $\kk[G]$ is almost Gorenstein if and only if one of the following holds:
    \begin{enumerate}
        \item $G$ is the complete bipartite graph $K_{2,m}$ for some $m\ge2$.
        \item $G$ is matching-covered and every $\{x,y\}\in\Fill(G)$, with $x\in X$ and $y\in Y$, satisfies
              \begin{equation*}
                  N_G(x)=Y\setminus\{y\},\qquad N_G(y)=X\setminus\{x\}.
              \end{equation*}
    \end{enumerate}
\end{thm}

Our second main result concerns the symmetry of the $h$-vector.
It applies to arbitrary bipartite graphs and shows that, in degree at least two, almost Gorensteinness forces symmetry except possibly between the linear and next-to-leading coefficients.

\begin{thm}[\Cref{cor:palindrome}]\label{thm:palindrome}
    Let $G$ be a bipartite graph such that $\kk[G]$ is almost Gorenstein, and write its $h$-polynomial as $h(\kk[G];t)=\sum_{i=0}^s h_it^i$ with $h_s\neq 0$.
    If $s\ge 2$,
    then one has
    \[
        h_i=h_{s-i}\qquad\left(i=0\text{ or }2\le i\le\left\lfloor s/2\right\rfloor\right).
    \]
\end{thm}

For a $2$-connected matching-covered bipartite graph $G$ with bipartition $X\sqcup Y$ and $|X|=|Y|=:n\ge3$, the $h$-polynomial has degree $s=n-1$.
Let $\widehat G$ be the graph obtained from $G$ by adding
all edges in $\Fill(G)$, and put $q:=|\Fill(G)|$.
In this case, the proof gives the more precise identity
\[
    h(\kk[\widehat G];t)=h(\kk[G];t)+qt
\]
whenever $\kk[G]$ is almost Gorenstein (see \eqref{eq:closure-h}).
Since $\kk[\widehat G]$ is Gorenstein (\cite[Proposition~5.3]{hatasa2026palindromicity}), this identity explains the coefficient symmetry in \Cref{thm:palindrome}.
This theorem gives an affirmative answer to \cite[Question~4.8]{hatasa2025pseudo} when $s\ge2$.

As applications of \Cref{thm:main}, we recover known results on almost Gorenstein edge rings for two families of graphs (see \Cref{ex:deleted-matching,ex:theta}).
Moreover, \Cref{ex:same-h} gives two bipartite graphs whose edge rings have the common $h$-vector $(1,5,10,6,1)$, although only one of the rings is almost Gorenstein.
Thus, almost Gorensteinness is not determined by the $h$-vector even among bipartite edge rings, answering \cite[Problem~4.6]{higashitani2016almost} within this class.

\;

This paper is organized as follows.
In \Cref{sec:preliminaries}, we recall basic facts on almost Gorenstein graded rings and bipartite edge rings.
In \Cref{sec:two-connected}, we prove our main theorems and present their applications.

\section{Preliminaries}\label{sec:preliminaries}

Throughout, $\kk$ is a field and all graphs are finite and simple.

\subsection{Almost Gorenstein graded rings}\label{subsec:almost-gorenstein}
We first recall some basic notions and invariants associated with standard graded Cohen--Macaulay domains (see, e.g., \cite{bruns1998cohen} for background).

Let $R=\bigoplus _{k\ge 0}R_k$ be a standard graded Cohen--Macaulay domain with $R_0=\kk$ and $d:=\dim R$.
Denote the unique homogeneous maximal ideal of $R$ by $\frakm:=\bigoplus_{k\ge 1} R_k$, and let $\omega_R$ be the canonical module of $R$.
For a finitely generated graded $R$-module $M$, write $\mu(M)$ for its minimal number of homogeneous generators and $e(M)$ for its multiplicity with respect to $\frakm$.
Note that the inequality $e(M)\ge \mu (M)$ holds if $M$ is Cohen--Macaulay.
Our convention for shifts is $M(b)_j:=M_{b+j}$.
Put
\[
    a(R):=-\min\{j:(\omega_R)_j\ne0\},\qquad
    \type(R):=\mu(\omega_R).
\]
The integers $a(R)$ and $\type(R)$ are called the \textit{$a$-invariant} and \textit{Cohen--Macaulay type} of $R$, respectively.

The \textit{Hilbert series} of $R$ is the formal power series $\Hilb_R(t):=\sum_{k\ge 0}(\dim_\kk R_k)t^k$.
It can be written in the form
\[
    \Hilb_R(t)=\frac{h_0+h_1t+\cdots+h_st^s}{(1-t)^d},
\]
where $h_s \neq 0$.
We call the polynomial $h_0+h_1t+\cdots+h_st^s$ the \textit{$h$-polynomial} of $R$, denote it by $h(R;t)$, and call the sequence $(h_0,h_1,\ldots,h_s)$ the \textit{$h$-vector} of $R$.
The index $s$ is called the \textit{socle degree} of $R$.
We have $s=d+a(R)$ and $h_s=\dim_\kk (\omega_R)_{-a(R)}$ (see \cite[Section 4.4]{bruns1998cohen}).
We call $R$ \textit{pseudo-Gorenstein} if $h_s=1$.

\begin{defi}[{\cite[Definition~1.5]{goto2015almost}}]\label[defi]{def:almost-gorenstein}
    We call $R$ \textit{almost Gorenstein} if there exists an exact sequence of graded $R$-modules
    \begin{equation}\label{eq:canonical-sequence}
        0\longrightarrow R\longrightarrow\omega_R(-a(R))\longrightarrow C\longrightarrow0
    \end{equation}
    such that $C=0$ or $\mu(C)=e(C)$.
    A nonzero Cohen--Macaulay graded module $C$ satisfying the equality $\mu(C)=e(C)$ is called an \textit{Ulrich module}.
\end{defi}
Under our assumptions on $R$, there always exists a degree-preserving injection from $R$ to $\omega_R(-a(R))$ (\cite[Proposition 2.2]{higashitani2016almost}), and its cokernel $C$ is Cohen--Macaulay of dimension $d-1$ if it is nonzero.
The following proposition gives a numerical criterion for almost Gorensteinness:
\begin{prop}[{\cite[Propositions~2.3, 2.4 and Corollary~2.7]{higashitani2016almost}}]\label[prop]{prop:canonical-cokernel}
    Take an exact sequence as in \eqref{eq:canonical-sequence} and let $(h_0,\ldots,h_s)$ be the $h$-vector of $R$.
    Then one has
    \begin{equation*}
        \Hilb_C(t)=\frac{\sum_{j=0}^{s-1}\left(\sum_{i=0}^j(h_{s-i}-h_i)\right)t^j}{(1-t)^{d-1}},
        \quad e(C)=\sum_{j=0}^{s-1}\sum_{i=0}^j(h_{s-i}-h_i),
        \quad \mu(C)=\type(R)-1.
    \end{equation*}
    Hence $R$ is almost Gorenstein if and only if $\sum_{j=0}^{s-1}\sum_{i=0}^j(h_{s-i}-h_i)=\type(R)-1$.
    Thus, this criterion is independent of the choice of the embedding $R\to \omega_R(-a(R))$ in \eqref{eq:canonical-sequence}.
\end{prop}

Since the $h$-vector and Cohen--Macaulay type of a standard graded normal affine semigroup ring depend only on its semigroup, \Cref{prop:canonical-cokernel} shows that almost Gorensteinness is independent of the coefficient field.
Thus, the following results, originally proved over an algebraically closed field of characteristic zero, hold over any field:
\begin{thm}[{\cite[Theorems~4.1 and~4.7]{higashitani2016almost}}]\label[thm]{thm:small-degree-facts}
    Let $R$ be a standard graded normal affine semigroup ring over $\kk$ with $h$-vector $(h_0,\ldots,h_s)$.
    Then the following hold:
    \begin{enumerate}
        \item If $s=1$, then $R$ is almost Gorenstein.
        \item If $R$ is almost Gorenstein and $s\ge2$, then $R$ is pseudo-Gorenstein, that is, $h_s=1$.
    \end{enumerate}
\end{thm}

The following lemma will be used in the proof of our main theorem:
\begin{lem}\label[lem]{lem:ulrich-multiplication}
    Let $R$ be a standard graded $\kk$-algebra, and let $M$ be a nonzero graded Ulrich $R$-module of dimension $\delta$.
    If $M_j=0$ for $j\le0$, then we have
    \[
        \dim_\kk (\frakm M)_2=\delta\dim_\kk M_1.
    \]
\end{lem}
\begin{proof}
    Extending scalars to $\kk(t)$ preserves the relevant dimensions, multiplicity, and number of generators, so we may assume that $\kk$ is infinite.
    Choose a linear system of parameters $\ell_1,\ldots,\ell_\delta$ for $M$ and put $Q=(\ell_1,\ldots,\ell_\delta)$.
    Since $M$ is Cohen--Macaulay and Ulrich,
    \[
        \dim_\kk M/QM =e(M)=\mu(M)=\dim_\kk M/\frakm M.
    \]
    The natural surjection $M/QM\longrightarrow M/\frakm M$ is an isomorphism, so $\frakm M=QM$.
    Moreover, $M$ is graded free over $\kk[\ell_1,\ldots,\ell_\delta]$.
    Since $M_j=0$ for $j\le0$, the degree-one elements of a homogeneous free basis of $M$ form a $\kk$-basis of $M_1$.
    Thus, we have
    \[
        \dim_\kk (\frakm M)_2= \dim_\kk(QM)_2=\delta\dim_\kk M_1.\qedhere
    \]
\end{proof}

\subsection{Bipartite edge rings}\label{subsec:edge-rings}

We next fix the graph-theoretic notation and recall the definition and some properties of edge rings.
For background on graph theory and edge rings, see \cite{diestel2017graph}, \cite[Section 10]{villarreal2001monomial}, and \cite[Section 5]{herzog2018binomial}.

Let $G$ be a simple connected graph on vertex set $V(G)$ with edge set $E(G)$.
For a vertex $v$, let $N_G(v)$ be its neighborhood and $d_G(v):=|N_G(v)|$.
For $T\subset V(G)$, write $N_G(T):=\bigcup_{v\in T}N_G(v)$, let $G[T]$ be the induced subgraph, and set $G\setminus T:=G[V(G)\setminus T]$.
We write $G\setminus v$ for $G\setminus\{v\}$.
We say that $G$ is \textit{$2$-connected} if $|V(G)|\ge 3$ and $G\setminus v$ is connected for any $v\in V(G)$.

In what follows, we assume that $G$ is bipartite with bipartition $V(G)=X\sqcup Y$.
We identify pairs in $X\times Y$ with unordered pairs of vertices.
For a subset $F\subset X\times Y$, let $G+F$ denote the graph on $V(G)$ with edge set $E(G)\cup F$.
For $f\in X\times Y$, we write $G+f$ for $G+\{f\}$.
For a nonempty set $T\subset X$, let $B(T)$ be the graph on $T\cup N_G(T)$ consisting of all edges incident with $T$.
We call $T$ \textit{acceptable} if $B(T)$ is connected and $G\setminus(T\cup N_G(T))$ is connected with at least one edge.
The same terminology applies to subsets of $Y$.

For an unordered pair $e=\{u,v\}$ of distinct vertices, put $\rho(e)=\eb_u+\eb_v\in\ZZ^{V(G)}$, where $\eb_i$ denotes the unit vector indexed by $i$.
Let $\kk[t_v : v\in V(G)]$ be the polynomial ring and let $\tb^{U}:=\prod_{u\in U}t_u$ for $U\subset V(G)$.
The \textit{edge ring} $\kk[G]$ and \textit{edge polytope} $P_G$ of $G$ are defined as
\[
    \kk[G]:=\kk[\tb^e:e\in E(G)]\subset \kk[t_v : v\in V(G)],\qquad
    P_G:=\conv\{\rho(e):e\in E(G)\}\subset \RR^{V(G)}.
\]
Each edge monomial $\tb^e$ has degree one.
The same definitions of the edge ring and edge polytope apply to graphs that are not bipartite.
Set
\[
    S_G=\ZZ_{\ge0}\{\rho(e):e\in E(G)\},\qquad
    \mathcal C_G=\RR_{\ge0}\{\rho(e):e\in E(G)\}.
\]
If $G$ is connected and bipartite with at least one edge, then $\kk[G]=\kk[S_G]$ is normal (\cite{ohsugi1998normal,simis1998integral}), hence Cohen--Macaulay (\cite{hochster1972rings}), and
\begin{equation}\label{eq:edge-lattice}
    \dim\kk[G]=|X|+|Y|-1,\qquad
    \ZZ S_G=\left\{z\in\ZZ^{X\sqcup Y}:\sum_{x\in X}z_x=\sum_{y\in Y}z_y\right\}.
\end{equation}
Normality gives $S_G=\mathcal C_G\cap\ZZ S_G$, so the degree-$k$ monomials of $\kk[G]$ correspond bijectively to the lattice points in $kP_G$.
Moreover, the canonical module of $\kk[G]$ is the monomial ideal with exponent vectors in $S_G\cap\relint(\mathcal C_G)$ (\cite[Theorem~6.3.5]{bruns1998cohen}).
Here, for a polytope or a cone $Q$, let $\relint(Q)$ denote its relative interior.

\begin{prop}[{\cite[Theorem~1.7]{ohsugi1998normal}}]\label[prop]{prop:facets}
    Let $G$ be a $2$-connected bipartite graph with bipartition $X\sqcup Y$.
    Every coordinate hyperplane supports a facet of $\mathcal C_G$.
    The other facet inequalities are
    \begin{equation*}
        L_T(z):=\sum_{y\in N_G(T)}z_y-\sum_{x\in T}z_x\ge0
    \end{equation*}
    for acceptable sets $T\subset X$.
    Consequently, $z\in\ZZ S_G$ belongs to $\relint(\mathcal C_G)$ if and only if all its coordinates and all its acceptable set forms $L_T(z)$ are positive.
\end{prop}
Note that the form $L_T$ is non-negative on the cone $\mathcal C_G$ for every $T\subset X$, even when $T$ is not acceptable.

A \textit{perfect matching} is a set of pairwise disjoint edges covering all vertices.
A graph is \textit{matching-covered} if it is connected, has at least one edge, and every edge belongs to a perfect matching.
The following theorem characterizes $2$-connected matching-covered bipartite graphs in terms of neighborhood sizes and pseudo-Gorensteinness.
\begin{thm}[{cf. \cite[Theorems~6.1.5 and 6.1.6]{asratian1998bipartite} and \cite[Theorem~3.6]{hatasa2025pseudo}}]\label{matching-covered}
    Let $G$ be a $2$-connected bipartite graph with bipartition $V(G)=X\sqcup Y$.
    Then the following are equivalent:
    \begin{enumerate}
        \item $G$ is matching-covered;
        \item $|X|=|Y|$ and $|N_G(T)|\ge |T|+1$ for every non-empty subset $T\subsetneq X$;
        \item $\kk[G]$ is pseudo-Gorenstein.
    \end{enumerate}
    In particular, if $G$ is matching-covered, then $a(\kk[G])=-|X|=-|Y|$ and
    the all-ones vector $\one:=(1,\ldots,1)\in \ZZ^{V(G)}$
    is the unique lattice point in $\relint(|X|P_G)$.
\end{thm}

We next recall a construction that produces a Gorenstein edge ring by adding certain missing edges to a matching-covered graph.
An acceptable set $T$ is \textit{tight} if $|N_G(T)|-|T|=1$.
\begin{defi}[{\cite[Section~2.1 and Proposition~5.3]{hatasa2026palindromicity}}]\label[defi]{def:fill}
    For a $2$-connected matching-covered bipartite graph $G$ with bipartition $X\sqcup Y$, define
    \[
        \Fill(G):=\left\{\{x,y\}\in X\times Y:\begin{array}{l}
            \{x,y\}\notin E(G),                                    \\
            \text{no tight acceptable }T\subset X\text{ satisfies} \\
            x\in T\text{ and }y\notin N_G(T)
        \end{array}\right\}.
    \]
    We call its elements \textit{filling edges} and put $\widehat G:=G+\Fill(G)$.
\end{defi}

\begin{prop}[{\cite[Propositions~3.3, 3.4 and~5.3]{hatasa2026palindromicity}}]\label[prop]{prop:fill-data}
    Let $G$ be a $2$-connected bipartite graph with bipartition $V(G)=X\sqcup Y$.
    Write $h(\kk[G];t)=\sum_{i=0}^s h_it^i$ with $h_s\ne0$.
    If $G$ is matching-covered,
    then
    \begin{enumerate}
        \item $\Fill(G)=\{f\in (X\times Y)\setminus E(G):\one+\rho(f)\in\relint(\mathcal C_G)\}$,
        \item $|\Fill(G)|=h_{s-1}-h_1$,
        \item $\widehat G$ is $2$-connected and matching-covered, $\kk[\widehat G]$ is Gorenstein, and $\deg h(\kk[\widehat G];t)=s$.
    \end{enumerate}
    In particular, $\kk[G]$ is Gorenstein if and only if $\Fill(G)=\emptyset$.
\end{prop}
We refer to $\widehat G$ as the \textit{Gorenstein closure} of $G$.

Let $G$ be a bipartite graph.
An edge $e\in E(G)$ is called a \textit{bridge} if removing $e$ increases the number of connected components of $G$.
A \textit{block} means a maximal $2$-connected subgraph or a bridge.
When $G$ has maximal $2$-connected components $B_1,\ldots,B_m$ and bridges $e_1,\ldots,e_\ell$, we have
\[
    \kk[G]\cong \left(\kk[B_1]\otimes_\kk \cdots \otimes_\kk \kk[B_m]\right)[x_1,\ldots,x_\ell],
\]
where $x_1,\ldots,x_\ell$ are indeterminates of degree one (see \cite[Proposition~10.1.48]{villarreal2001monomial}).
Since each $B_i$ contains a cycle, $\kk[B_i]$ is not a polynomial ring, and hence $e(\kk[B_i])>1$.
Thus, the following holds by \cite[Corollaries~1.5 and 3.10]{matsushita2024multiplicity}:
\begin{prop}\label[prop]{prop:2-connected}
    With the above notation, $\kk[G]$ is Gorenstein if and only if each $\kk[B_i]$ is Gorenstein.
    Moreover, if $\kk[G]$ is not Gorenstein but almost Gorenstein, then $m=1$ and $\kk[B_1]$ is almost Gorenstein.
\end{prop}
Therefore, it remains to characterize almost Gorenstein edge rings of $2$-connected bipartite graphs.

\section{Characterization and consequences}\label{sec:two-connected}

\subsection{Key lemmas}\label{subsec:lemmas}
We first establish two lemmas used in the proof of our main theorems.
Until the proof of \Cref{thm:main}, $G$ is a $2$-connected matching-covered bipartite graph with $|X|=|Y|=:n$ and $n\ge3$.
Set $R:=\kk[G]$, $\frakm:=\bigoplus_{k\ge 1} R_k$, $q:=|\Fill(G)|$, and $\tau:=\tb^{V(G)}=\prod_{v\in V(G)}t_v$.
Let $(h_0,\ldots,h_s)$ be the $h$-vector of $R$.
Consider the exact sequence
\begin{equation}\label{eq:monomial-sequence}
    0\longrightarrow R\xrightarrow{\,1\mapsto\tau\,}\omega_R(n)\longrightarrow C\longrightarrow0.
\end{equation}
Note that we have $C_j=0$ for $j<0$.
Moreover, since $h_0=h_s=1$ by \Cref{matching-covered}, the Hilbert-series formula in
\Cref{prop:canonical-cokernel} yields
$C_0=0$ and $\dim_\kk C_1=h_{s-1}-h_1$.

By \Cref{prop:fill-data} (i), the classes
$[\tau\tb^f]$ for $f\in\Fill(G)$
are nonzero and linearly independent in $C_1$.
Since their number is $q=h_{s-1}-h_1=\dim_\kk C_1$
by \Cref{prop:fill-data} (ii),
they form a $\kk$-basis of $C_1$.
Consequently, $(\frakm C)_2$ is spanned by the classes
$[\tau\tb^f\tb^e]$ for $f\in\Fill(G)$ and $e\in E(G)$.
The following lemma examines which of these classes
are nonzero and distinct, and gives a lower bound
for $\dim_\kk(\frakm C)_2$:

\begin{lem}\label[lem]{lem:count}
    For $f=\{x,y\}\in (X\times Y)\setminus E(G)$, define
    \[
        P_f:=N_G(y),\quad Q_f:=N_G(x),\quad
        U_f:=X\setminus(P_f\cup\{x\}),\quad
        V_f:=Y\setminus(Q_f\cup\{y\}).
    \]
    Then, we have
    \begin{equation}\label{eq:w-bound}
        \dim_\kk (\frakm C)_2 \ge(2n-2)q+
        \frac{1}{2}\sum_{f\in\Fill(G)}\left(
        \sum_{u\in U_f}(d_G(u)-2)+\sum_{v\in V_f}(d_G(v)-2)\right).
    \end{equation}
    Moreover, if $\dim_\kk (\frakm C)_2 =(2n-2)q$, then $U_f=V_f=\emptyset$ for every $f\in\Fill(G)$.
\end{lem}
\begin{proof}
    For $f=\{x,y\}\in\Fill(G)$ and $e=\{u,v\}\in E(G)$, with $u\in X$ and $v\in Y$, consider $[\tau\tb^f\tb^e]\in C_2$.
    This class vanishes precisely when $\tb^f\tb^e\in R_2$.
    The exponent vector of $\tb^f\tb^e$ is $\rho(f)+\rho(e)$, which has coordinate sum two on each side.
    There are only two possible pairings of the vertices $x,u$ with $y,v$, counted with multiplicity.
    Since $f$ is not an edge, the class vanishes exactly when $\{u,y\}$ and $\{x,v\}$ are both edges, that is, when $u\in P_f$ and $v\in Q_f$.

    A nonzero class has at most two representations as $[\tau\tb^g\tb^{e'}]$ with $g\in\Fill(G)$ and $e'\in E(G)$.
    Besides $(f,e)$, the only possible representation uses one of the cross pairs $\{u,y\}$, $\{x,v\}$ as a filling edge and the other as an edge of $G$.
    If $e$ meets $f$, there is only one representation.
    Define
    \[
        \lambda_{f,e}=\begin{cases}
            0,   & [\tau\tb^f\tb^e]=0,                                           \\
            1,   & [\tau\tb^f\tb^e]\ne0\text{ and its representation is unique}, \\
            1/2, & [\tau\tb^f\tb^e]\ne0\text{ and it has two representations}.
        \end{cases}
    \]
    Nonzero monomial classes with distinct exponent vectors are linearly independent.
    The classes $[\tau\tb^f\tb^e]$ span $(\frakm C)_2$, so
    \begin{equation*}
        \dim_\kk (\frakm C)_2 =\sum_{f\in\Fill(G)}\sum_{e\in E(G)}\lambda_{f,e}.
    \end{equation*}

    For fixed $f$, put
    \[
        B_f:=E(G)\cap\bigl((U_f\times Q_f)\cup(P_f\times V_f)\bigr),
    \]
    where pairs are identified with their unordered edges.
    The possible weights according to the location of $e$ are
    \[
        \begin{array}{c|c}
            \text{location of }e               & \lambda_{f,e}    \\\hline
            \text{incident with }x\text{ or }y & 1                \\
            \text{between }P_f\text{ and }Q_f  & 0                \\
            \text{between }U_f\text{ and }V_f  & 1                \\
            e\in B_f                           & 1/2\text{ or }1.
        \end{array}
    \]
    The degree sum over $U_f\cup V_f$ counts each edge between $U_f$ and $V_f$ twice and each edge of $B_f$ once.
    Consequently,
    \begin{align*}
        \sum_{e\in E(G)}\lambda_{f,e}
         & \ge d_G(x)+d_G(y)+\frac{1}{2}\left(
        \sum_{u\in U_f}d_G(u)+\sum_{v\in V_f}d_G(v)\right) \\
         & =2n-2+\frac{1}{2}\left(
        \sum_{u\in U_f}(d_G(u)-2)+\sum_{v\in V_f}(d_G(v)-2)\right).
    \end{align*}
    The last equality uses $d_G(x)=n-1-|V_f|$ and $d_G(y)=n-1-|U_f|$.
    Summing over $f$ proves \eqref{eq:w-bound}.

    Since $G$ is $2$-connected, every vertex has degree at least two.
    Moreover, every $f=\{x,y\}\in\Fill(G)$ satisfies
    \begin{equation}\label{eq:fill-endpoint-degree}
        d_G(x)\ge3,\qquad d_G(y)\ge3.
    \end{equation}
    To see this, let $f=\{x,y\}\in\Fill(G)$.
    The non-negative form $L_{\{x\}}$ on $\mathcal C_G$ is not identically zero, since connectedness gives an edge between $X\setminus\{x\}$ and $N_G(x)$.
    Since $\one+\rho(f)\in\relint(\mathcal C_G)$, we have
    \[
        0<L_{\{x\}}(\one+\rho(f))=d_G(x)-2.
    \]
    The symmetric argument applies to $y$.

    If the equality $\dim_\kk (\frakm C)_2 =(2n-2)q$ holds, all terms $d_G(u)-2$ and $d_G(v)-2$ above vanish, and every edge of $B_f$ has weight $1/2$.
    Suppose $U_f\cup V_f\ne\emptyset$.
    Connectedness gives an edge between $U_f\cup V_f$ and its complement, so $B_f\ne\emptyset$.
    If $e=\{u,v\}\in B_f$ with $u\in U_f$ and $v\in Q_f$, its weight $1/2$ forces $\{u,y\}\in\Fill(G)$, since $\{x,v\}\in E(G)$.
    But an endpoint of a filling edge has degree at least three, whereas $d_G(u)=2$.
    For $u\in P_f$ and $v\in V_f$, the same argument uses $\{x,v\}\in\Fill(G)$ and contradicts $d_G(v)=2$.
    Thus $U_f=V_f=\emptyset$.
\end{proof}

We say that $f=(x,y)\in (X\times Y)\setminus E(G)$ satisfies the \textit{neighborhood condition} if $U_f=V_f=\emptyset$, equivalently, $N_G(x)=Y\setminus \{y\}$ and $N_G(y)=X\setminus \{x\}$.
The following lemma describes how the Hilbert series changes when a missing edge satisfying the neighborhood condition is added to $G$:

\begin{lem}\label[lem]{lem:single-edge}
    Let $G$ be a connected bipartite graph with $|X|=|Y|=:n$.
    Suppose that $f=(x,y)\in (X\times Y)\setminus E(G)$ satisfies the neighborhood condition.
    Then
    \[
        \Hilb_{\kk[G+f]}(t)-\Hilb_{\kk[G]}(t)=\frac{t}{(1-t)^{2n-1}},
    \]
    and hence $h(\kk[G+f];t)=h(\kk[G];t)+t$.
\end{lem}
\begin{proof}
    Whenever a product of edge monomials of $G+f$ contains both $\tb^f$ and $\tb^{\{u,v\}}$, with $u\in X\setminus\{x\}$ and $v\in Y\setminus\{y\}$, we can replace their product using
    \[
        (t_xt_y)(t_ut_v)=(t_ut_y)(t_xt_v).
    \]
    The two monomials on the right correspond to edges of $G$ by the neighborhood condition.
    Each replacement reduces the number of factors $\tb^f$ by one.
    For a monomial outside $\kk[G]$, at least one such factor must remain.
    Repeating the replacement therefore expresses every such monomial in the form
    \begin{equation}\label{eq:new-monomials}
        (t_xt_y)^r\prod_{u\in X\setminus\{x\}}(t_ut_y)^{a_u}
        \prod_{v\in Y\setminus\{y\}}(t_xt_v)^{b_v},
        \qquad r\ge1,\quad a_u,b_v\ge0.
    \end{equation}
    The exponents $a_u,b_v$, and then $r$, are uniquely determined by the exponent vector.
    Conversely, for a monomial of degree $k$ in \eqref{eq:new-monomials}, the sum of its $x$- and $y$-exponents is $k+r>k$.
    Every degree-$k$ monomial of $\kk[G]$ has that sum at most $k$, since no edge of $G$ meets both $x$ and $y$.
    Thus all monomials in \eqref{eq:new-monomials} lie outside $\kk[G]$.
    Their generating series is
    \[
        (t+t^2+\cdots)(1+t+t^2+\cdots)^{2n-2}
        =\frac{t}{(1-t)^{2n-1}}.\qedhere
    \]
\end{proof}

\subsection{Proofs of the main results}\label{subsec:proofs}

We now prove \Cref{thm:main}.
The proof also gives the numerical formulas recorded in \Cref{cor:palindrome}.

\begin{proof}[Proof of \Cref{thm:main}]
    Put $R:=\kk[G]$ and let $s$ be its socle degree.
    Since $G$ is $2$-connected, it has a cycle and $s\ge1$.
    The case $s=1$ occurs exactly when $G=K_{2,m}$ for some $m\ge2$; this follows from \cite[Lemma~4.2 and Theorem~4.6 and its proof]{ohsugi1999quadratic}, since $G$ is $2$-connected.
    By \Cref{thm:small-degree-facts} (i), these rings are almost Gorenstein.
    In what follows, we assume $s\ge2$.

    Assume first that $R$ is almost Gorenstein.
    By \Cref{thm:small-degree-facts} (ii), we have $h_s=1$, so $G$ is matching-covered by \Cref{matching-covered}.
    Put $n:=|X|=|Y|$ and $q:=|\Fill(G)|$.
    Then \eqref{eq:edge-lattice} and \Cref{matching-covered} give $\dim R=2n-1$, $a(R)=-n$, and $s=n-1$.
    In particular, $n=s+1\ge3$.

    If $\Fill(G)=\emptyset$, equivalently, $\kk[G]$ is Gorenstein, the neighborhood condition holds automatically.
    Otherwise the cokernel $C$ in \eqref{eq:monomial-sequence} is nonzero and Ulrich by \Cref{prop:canonical-cokernel}, and \Cref{lem:ulrich-multiplication} gives $\dim_\kk (\frakm C)_2 =(\dim \kk[G]-1)\cdot \dim_\kk C_1=(2n-2)q$.
    By \Cref{lem:count}, $U_f=V_f=\emptyset$ for every $f\in\Fill(G)$.

    Conversely, suppose that (ii) holds, that is, $G$ is matching-covered and every $f\in \Fill(G)$ satisfies the neighborhood condition.
    Two different edges of $\Fill(G)$ cannot share an endpoint, since each such endpoint has exactly one non-neighbor across the bipartition.
    Thus $\Fill(G)$ is a matching in $\widehat G$.
    Adjoin its edges one at a time.
    For every remaining filling edge its endpoint neighborhoods are unchanged, so \Cref{lem:single-edge} applies at each step and yields
    \begin{equation}\label{eq:closure-h}
        h(\kk[\widehat G];t)=h(R;t)+qt.
    \end{equation}
    The ring $\kk[\widehat G]$ is Gorenstein by \Cref{prop:fill-data} (iii), so its $h$-polynomial is palindromic of degree $s=n-1$.
    If $q=0$, then $R$ is Gorenstein and hence almost Gorenstein.
    If $q>0$, \eqref{eq:fill-endpoint-degree} and the neighborhood condition give $n\ge4$.
    The palindromicity of $h(\kk[\widehat G];t)=h(R;t)+qt$ and \Cref{prop:canonical-cokernel} imply
    \begin{equation*}
        \Hilb_C(t)=\frac{q(t+t^2+\cdots+t^{n-3})}{(1-t)^{2n-2}},
        \qquad e(C)=(n-3)q.
    \end{equation*}
    For each $f=\{x,y\}\in\Fill(G)$ and $1\le j\le n-3$, the exponent vector
    \[
        \one+j\rho(f)
    \]
    lies in $\relint(\mathcal C_G)$.
    Indeed, its coordinates are positive.
    For an acceptable $T\subset X$, the only way that the added coordinates decrease $L_T$ is $x\in T$ and $y\notin N_G(T)$.
    The neighborhood condition then forces $T=\{x\}$, and
    \[
        L_{\{x\}}(\one+j\rho(f))=n-2-j>0.
    \]
    For every other acceptable set, the value is at least $|N_G(T)|-|T|\ge1$.
    Normality shows that these lattice vectors give canonical monomials.

    All these monomials are minimal generators of $\omega_R$.
    If subtracting an edge vector left an interior vector, positivity of every coordinate would force both endpoints of the subtracted edge to lie in $\{x,y\}$.
    But $f$ is not an edge of $G$.
    The monomial $\tau$ is also a minimal generator of $\omega_R$.
    These monomials have distinct exponent vectors, so
    \[
        \type(R)\ge1+(n-3)q.
    \]
    On the other hand, \Cref{prop:canonical-cokernel} gives
    \[
        \type(R)-1=\mu(C)\le e(C)=(n-3)q.
    \]
    Equality follows, proving almost Gorensteinness.
\end{proof}

The proof of \Cref{thm:main} also yields our second main result:

\begin{cor}\label[cor]{cor:palindrome}
    Let $G$ be a bipartite graph such that $\kk[G]$ is
    almost Gorenstein, and write
    $h(\kk[G];t)=\sum_{i=0}^s h_it^i$ with $s\ge2$.
    Then $q:=h_{s-1}-h_1\ge0$ and
    $h(\kk[G];t)+qt$ is palindromic.
    Moreover, $\type(\kk[G])=1+(s-2)q$.
\end{cor}
\begin{proof}
    If $\kk[G]$ is Gorenstein, the assertions are clear.
    Otherwise, by \Cref{prop:2-connected}, $\kk[G]$
    is a polynomial extension of $\kk[B]$ for a
    $2$-connected bipartite graph $B$.
    Since $s\ge2$, \Cref{thm:main} implies that
    $B$ is matching-covered.
    The assertions now follow from the proof of
    \Cref{thm:main} and \Cref{prop:fill-data} (ii).
\end{proof}

This gives the coefficient symmetry asked for in \cite[Question~4.8]{hatasa2025pseudo} when $s\ge2$.
The degree restriction is necessary: $\kk[K_{2,m}]$ is almost Gorenstein with $h$-polynomial $1+(m-1)t$, which has $h_0\ne h_s$ when $m\ge3$.

\subsection{Examples}\label{subsec:examples}

We first recover the almost Gorenstein criteria for two known families.

\begin{ex}[{\cite{hatasa2025pseudo}}]\label[ex]{ex:deleted-matching}
    Let $m,n\ge3$, and let $G$ be obtained from $K_{m,n}$ by deleting a matching.
    Then $\kk[G]$ is almost Gorenstein if and only if $m=n$.
    Indeed, $G$ is $2$-connected and is not a graph of the form $K_{2,r}$, so \Cref{thm:main} makes $m=n$ necessary.
    If $m=n$, then $G$ is matching-covered by \Cref{matching-covered}.
    Every missing edge satisfies the neighborhood condition, and \Cref{thm:main} gives the assertion.
\end{ex}

The next example uses a bipartite realization of an edge ring of a non-bipartite graph.

\begin{ex}[{\cite{higashitani2023h,bhaskara2026h}}]\label[ex]{ex:theta}
    For $\ell\ge1$ and $k_1,\ldots,k_\ell\ge1$, let $\mathcal O(k_1,\ldots,k_\ell)$ consist of cycles of lengths $2k_1+1,\ldots,2k_\ell+1$ sharing a vertex $o$ and otherwise disjoint.
    Let $\Theta(k_1,\ldots,k_\ell)$ consist of internally disjoint $u$--$v$ paths of the same lengths (see \Cref{fig:odd-cycle-theta}).
    \begin{figure}[htbp]
        \centering
        \begin{tikzpicture}[
                x=1cm,y=1cm,line width=0.85pt,
                vertex/.style={circle,draw=black,fill=white,inner sep=0pt,
                        minimum size=4.2pt,line width=0.55pt}]
            \begin{scope}
                \coordinate (o) at (0,0);
                \foreach \i/\ang in {1/0, 2/110, 3/250}{
                        \begin{scope}[rotate=\ang]
                            \coordinate (a\i-1) at (-0.8,0.8);
                            \coordinate (a\i-2) at (-0.6,1.8);
                            \coordinate (a\i-3) at (0.6,1.8);
                            \coordinate (a\i-4) at (0.8,0.8);
                            \draw[dotted] (-0.2,1.8)--(0.2,1.8);
                            \draw (0,0)--(a\i-1)--(a\i-2)--(-0.3,1.8);
                            \draw (0.3,1.8)--(a\i-3)--(a\i-4)--(0,0);
                            \foreach \p in {a\i-1,a\i-2,a\i-3,a\i-4}
                            \node[vertex] at (\p) {};
                        \end{scope}
                    }
                \node[vertex,label=above:{$o$}] at (o) {};
                \draw[dotted] (255:1.5cm) arc[start angle=255, end angle=285, radius=1.5cm];
                \node at (0,-2.5) {$\mathcal O(k_1,\ldots,k_\ell)$};
                \node at (0,1.2) {\small $C_{2k_1+1}$};
                \node at (1,-0.4) {\small $C_{2k_\ell+1}$};
                \node at (-1,-0.4) {\small $C_{2k_2+1}$};
            \end{scope}
            \begin{scope}[xshift=5.8cm]
                \coordinate (u) at (0,0);
                \coordinate (v) at (4,0);
                \coordinate (a1) at (1,1.5);
                \coordinate (a2) at (3,1.5);
                \coordinate (b1) at (1,0.5);
                \coordinate (b2) at (3,0.5);
                \coordinate (c1) at (1,-1.5);
                \coordinate (c2) at (3,-1.5);
                \coordinate (d1) at (1,-0.5);
                \coordinate (d2) at (3,-0.5);
                \draw (u)--(a1)--(1.5,1.5); \draw (2.5,1.5)--(a2)--(v);
                \draw (u)--(b1)--(1.5,0.5); \draw (2.5,0.5)--(b2)--(v);
                \draw (u)--(c1)--(1.5,-1.5); \draw (2.5,-1.5)--(c2)--(v);
                \draw (u)--(d1)--(1.5,-0.5); \draw (2.5,-0.5)--(d2)--(v);
                \draw[dotted] (1.7,1.5)--(2.3,1.5); \draw[dotted] (1.7,0.5)--(2.3,0.5);
                \draw[dotted] (1.7,-1.5)--(2.3,-1.5); \draw[dotted] (1.7,-0.5)--(2.3,-0.5);
                \draw[dotted] (2,0.25)--(2,-0.25);
                \foreach \p in {a1,a2,b1,b2,c1,c2,d1,d2}
                \node[vertex] at (\p) {};
                \node[vertex,label=left:{$u$}] at (u) {};
                \node[vertex,label=right:{$v$}] at (v) {};
                \node at (2,-2.5) {$\Theta(k_1,\ldots,k_\ell)$};
                \draw [line width=0.015cm, decorate,decoration={brace,amplitude=5pt}](1,1.6) -- (3,1.6) node[midway,yshift=0.4cm] {\small $2k_1$};
                \draw [line width=0.015cm, decorate,decoration={brace,amplitude=5pt}](1,0.6) -- (3,0.6) node[midway,yshift=0.4cm] {\small $2k_2$};
                \draw [line width=0.015cm, decorate,decoration={brace,amplitude=5pt,mirror}](1,-0.6) -- (3,-0.6) node[midway,yshift=-0.4cm] {\small $2k_{\ell-1}$};
                \draw [line width=0.015cm, decorate,decoration={brace,amplitude=5pt,mirror}](1,-1.6) -- (3,-1.6) node[midway,yshift=-0.4cm] {\small $2k_{\ell}$};
            \end{scope}
        \end{tikzpicture}
        \caption{The graphs $\mathcal O(k_1,\ldots,k_\ell)$ and $\Theta(k_1,\ldots,k_\ell)$, where $C_{2k_i+1}$ in $\mathcal O(k_1,\ldots,k_\ell)$ denotes a cycle of length $2k_i+1$, while $2k_i$ in $\Theta(k_1,\ldots,k_\ell)$ denotes the number of internal vertices of the corresponding path.}
        \label{fig:odd-cycle-theta}
    \end{figure}
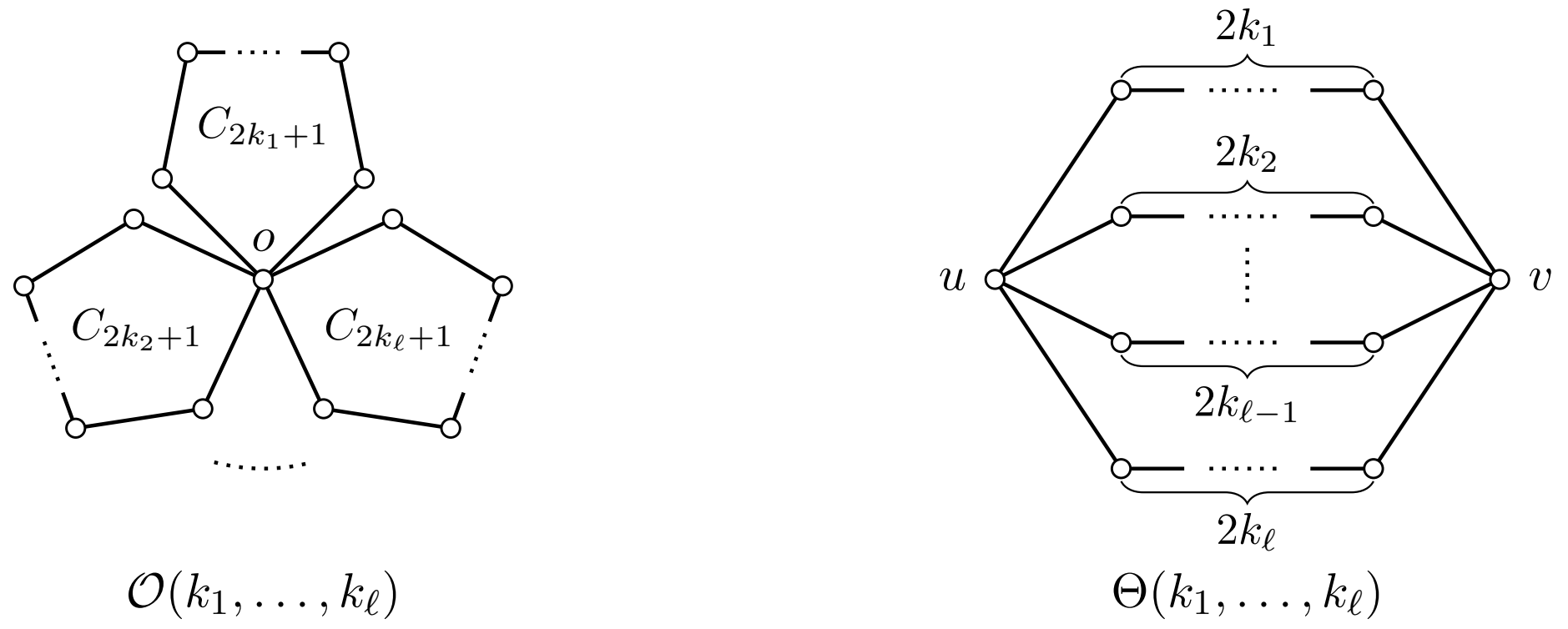

    Then we have $\kk[\mathcal O(k_1,\ldots,k_\ell)]\cong
        \kk[\Theta(k_1,\ldots,k_\ell)]$.
    Indeed, identifying $u$ and $v$ in $\Theta$ yields $\mathcal O$.
    Deleting the $u,v$-coordinates from $P_\Theta$ and the $o$-coordinate from $P_{\mathcal O}$ gives the same polytope.
    Both projections are lattice isomorphisms, since the deleted
    coordinates can be recovered using formulas with integer coefficients from the bipartition
    equations for $\Theta$ and the coordinate-sum equation for
    $\mathcal O$.

    For $\ell\le2$, $\Theta$ is a path or an even cycle, so its edge ring is Gorenstein.
    Suppose $\ell\ge3$.
    Then $\Theta$ is $2$-connected and matching-covered.
    Since every vertex other than $u$ and $v$ has degree two, \eqref{eq:fill-endpoint-degree} implies
    $\Fill(\Theta)\subseteq\{\{u,v\}\}$.
    It is straightforward to check that $\{u,v\}\in\Fill(\Theta)$, and hence $\Fill(\Theta)=\{\{u,v\}\}$.
    The pair $\{u,v\}$ satisfies the neighborhood condition if and only if $k_1=\cdots=k_\ell=1$.
    Thus, by \Cref{thm:main}, $\kk[\Theta]$ is almost Gorenstein if and only if $\ell\le2$ or $k_1=\cdots=k_\ell=1$.

\end{ex}

The last example shows that almost Gorensteinness is not determined by the $h$-vector, even among edge rings of $2$-connected matching-covered bipartite graphs.

\begin{ex}\label[ex]{ex:same-h}
    Let $X=\{x_1,\ldots,x_5\}$ and $Y=\{y_1,\ldots,y_5\}$.
    Define bipartite graphs $G$ and $G'$ with bipartition $X\sqcup Y$ by the following neighborhoods; each entry lists the indices of the $y$-vertices.
    \[
        \begin{array}{c|ccccc}
               & x_1         & x_2       & x_3       & x_4       & x_5     \\\hline
            G  & \{2,3,4,5\} & \{1,2,3\} & \{1,2,3\} & \{1,4\}   & \{1,5\} \\
            G' & \{3,4,5\}   & \{1,2,3\} & \{1,2,3\} & \{1,4,5\} & \{1,5\}
        \end{array}
    \]
    Thus $G'$ is obtained from $G$ by deleting $\{x_1,y_2\}$ and adding $\{x_4,y_5\}$.
    The two graphs are drawn in \Cref{fig:same-h}.
    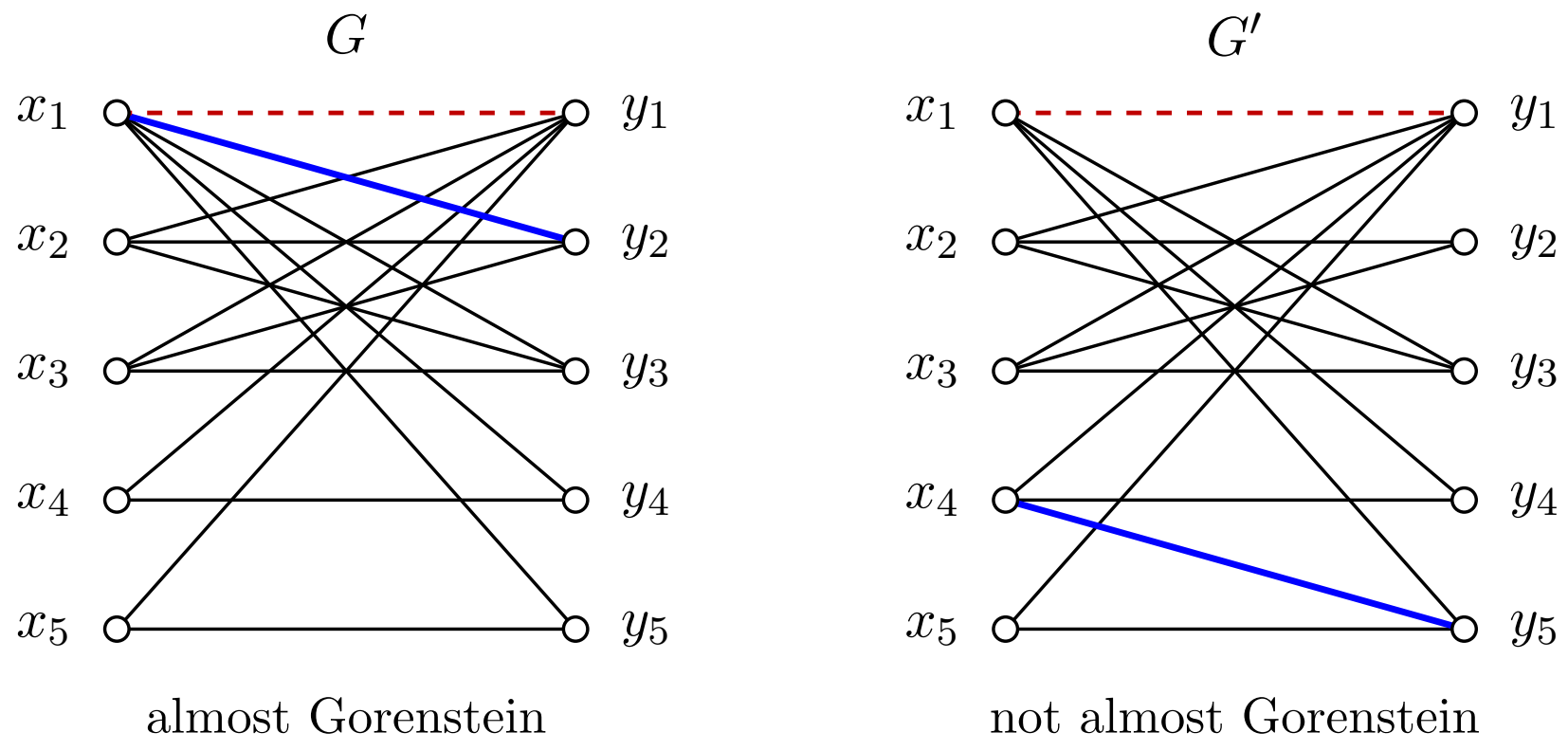
\begin{figure}[htbp]
        \centering
        \begin{tikzpicture}[x=1cm,y=1cm,line width=0.65pt]
            \begin{scope}
                \foreach \i in {1,...,5} {
                        \coordinate (x\i) at (0,{(5-\i)*0.9});
                        \coordinate (y\i) at (3.2,{(5-\i)*0.9});
                    }
                \foreach \i/\j in {1/2,1/3,1/4,1/5,2/1,2/2,2/3,3/1,3/2,3/3,4/1,4/4,5/1,5/5}
                \draw (x\i)--(y\j);
                \draw[blue,line width=1.3pt] (x1)--(y2);
                \draw[red!75!black,dashed,line width=0.9pt] (x1)--(y1);
                \foreach \i in {1,...,5} {
                        \draw[fill=white] (x\i) circle (2.4pt) node[left=5pt] {$x_\i$};
                        \draw[fill=white] (y\i) circle (2.4pt) node[right=5pt] {$y_\i$};
                    }
                \node at (1.6,4.15) {$G$};
                \node[font=\small] at (1.6,-0.6) {almost Gorenstein};
            \end{scope}
            \begin{scope}[xshift=6.2cm]
                \foreach \i in {1,...,5} {
                        \coordinate (x\i) at (0,{(5-\i)*0.9});
                        \coordinate (y\i) at (3.2,{(5-\i)*0.9});
                    }
                \foreach \i/\j in {1/3,1/4,1/5,2/1,2/2,2/3,3/1,3/2,3/3,4/1,4/4,4/5,5/1,5/5}
                \draw (x\i)--(y\j);
                \draw[blue,line width=1.3pt] (x4)--(y5);
                \draw[red!75!black,dashed,line width=0.9pt] (x1)--(y1);
                \foreach \i in {1,...,5} {
                        \draw[fill=white] (x\i) circle (2.4pt) node[left=5pt] {$x_\i$};
                        \draw[fill=white] (y\i) circle (2.4pt) node[right=5pt] {$y_\i$};
                    }
                \node at (1.6,4.15) {$G'$};
                \node[font=\small] at (1.6,-0.6) {not almost Gorenstein};
            \end{scope}
        \end{tikzpicture}
        \caption{The graphs $G$ and $G'$. Thick solid edges mark $\{x_1,y_2\}$ in $G$ and $\{x_4,y_5\}$ in $G'$. Dashed segments indicate the common filling pair $\{x_1,y_1\}$, which is not an edge of either graph. Crossings are not vertices.}
        \label{fig:same-h}
    \end{figure}

    Using \texttt{Macaulay2} (\cite{Macaulay2}), we obtain
    \[
        h(\kk[G];t)=h(\kk[G'];t)=1+5t+10t^2+6t^3+t^4.
    \]

    Both graphs are $2$-connected.
    Since $h_4=1$, both graphs are matching-covered by \Cref{matching-covered}.
    By \Cref{prop:fill-data} (ii), $|\Fill(G)|=|\Fill(G')|=h_3-h_1=1$.
    Since every vertex of $X\setminus\{x_1\}$ is adjacent to $y_1$ and $\{x_1\}$ is not tight in either graph, \Cref{def:fill} shows that $\{x_1,y_1\}$ belongs to both filling sets.
    Hence
    \[
        \Fill(G)=\Fill(G')=\{\{x_1,y_1\}\}.
    \]

    For $G$, the neighborhoods of $x_1$ and $y_1$ are $Y\setminus\{y_1\}$ and $X\setminus\{x_1\}$, respectively.
    For $G'$, $N_{G'}(x_1)=\{y_3,y_4,y_5\}$, so the required neighborhood equality fails.
    Thus, by \Cref{thm:main}, $\kk[G]$ is almost Gorenstein but $\kk[G']$ is not.
\end{ex}

The two rings in \Cref{ex:same-h} are normal Cohen--Macaulay homogeneous domains.
They therefore answer \cite[Problem~4.6]{higashitani2016almost} even when both domains are required to be bipartite edge rings.
Since their common $h$-polynomial becomes palindromic after adding $t$, the palindromicity condition in \Cref{cor:palindrome} does not imply almost Gorensteinness.

\subsection*{Acknowledgement}
The second author was partially supported by Grant-in-Aid for JSPS Fellows Grant JP25KJ0047.
The authors used OpenAI's ChatGPT and Codex to assist with drafting and editing the manuscript and to discuss some mathematical arguments.
The authors take full responsibility for the content of this paper.

\bibliographystyle{plain}
\bibliography{References}

\end{document}